\documentclass[12pt,a4paper]{article}
\usepackage{amsfonts}
\usepackage{bbm}
\usepackage{mathrsfs}

\usepackage{amssymb}

\usepackage{amsmath}

\newtheorem{theorem}{Theorem}[section]
\newtheorem{definition}[theorem]{Definition}
\newtheorem{lemma}[theorem]{Lemma}
\newtheorem{corollary}[theorem]{Corollary}

\begin{document}
\title{Composition-Diamond Lemma for Differential Algebras\footnote{Supported by the
NNSF of China (No.10771077) and the NSF of Guangdong Province
(No.06025062).} }
\author{
Yuqun Chen, Yongshan Chen  and Yu Li   \\
\\
{\small \ School of Mathematical Sciences}\\
{\small \ South China Normal University}\\
{\small \ Guangzhou 510631}\\
{\small \ P. R. China}\\
{\small \ yqchen@scnu.edu.cn} \\
{\small \ jackalshan@126.com} \\
{\small \ liyu820615@126.com}}
\date{}
\maketitle \noindent\textbf{Abstract:} In this paper, we establish
the Composition-Diamond lemma for free differential algebras. As
applications, we give Gr\"{o}bner-Shirshov bases for free
Lie-differential algebras and free commutative-differential
algebras, respectively.


\noindent {\bf AMS} Mathematics Subject Classification(2000): 16S15,
13P10, 16S32

\section{Introduction}

In 1962, A. I. Shirshov \cite{s62} invented a theory that is now
called the Gr\"{o}bner-Shirshov bases theory for Lie algebras
$Lie(X|S)$ presented by generators and defining relations. The main
technical notion of Shirshov's theory was a notion of composition
$(f,g)_w$ of two Lie polynomials, $f,g\in Lie(X)$ relative to some
associative word $w$. Based on it, he defined an infinite algorithm
of adding to some set $S$ of Lie polynomials of all non-trivial
compositions until one will get a set $S^c$ that is closed under
compositions, i.e., any non-trivial composition of two polynomials
from $S^c$ belongs to $S^c$ (and leading associative words $\bar{s}$
of polynomials $s\in S^c$ do not contain each other as subwords). In
addition, $S$ and $S^c$ generate the same ideal, i.e.,
$Id(S)=Id(S^c)$. $S^c$ is now called a Gr\"{o}bner-Shirshov basis of
$Lie(X|S)$. Then he proved the following lemma.

\emph{Let $Lie(X)\subset \mathbbm{k}\langle X\rangle$  be a free Lie
algebra over a field $\mathbbm{k}$ viewed as the algebra of Lie
polynomials in the free algebra $\mathbbm{k}\langle X\rangle$, and
$S$ a subset in $Lie(X)$. If $f\in Id(S)$, then $\bar{f}=u\bar{s}v$,
where $s\in S^c, \ u,v\in X^*, \ \bar{f},\bar{s}$ are leading
associative words of Lie polynomials $f,s$ correspondingly, and
$X^*$ the free monoid generated by $X$.}

He used the following easy corollary of his lemma.

\emph{$Irr(S^c)=\{[u] \ | \  u\neq{a\bar{s}b}, \ s\in{S^c},\
a,b\in{X^*}\}$ is a linear basis of the algebra
$Lie(X|S)=Lie(X)/Id(S)$, where $u$ is an associative Lyndon-Shirshov
word in $X^*$ and $[u]$ the corresponding non-associative
Lyndon-Shirshov word under Lie brackets $[xy]=xy-yx$.}

To define the Lie composition $(f,g)_w$ of two, say, monic Lie
polynomials, where $\bar{f}=ac, \ \bar{g}=cb, \ c\neq 1, \ a,b,c$
are associative words, and $w=acb$, A. I. Shirshov defines first the
associative composition $fb-ag$. Then he puts on $fb$ and $ag$
special brackets $[fb],[ag]$ in accordance with his paper
\cite{s58}. The result is $(f,g)_w=[fc]-[cg]$. Following \cite{s62},
one can easily get the same lemma for a free associative algebra:
Let $S\subset \mathbbm{k}\langle X\rangle$ and $S^c$ be as before.
If $f\in Id(S)$, then $\bar{f}=a\bar{s}b$ for some $s\in S^c, \
a,b\in X^*$. It was formulated later by L. A. Bokut \cite{b76} as an
analogue of Shirshov's Lie composition lemma, and by G. Bergman
\cite{b} under the name ``Diamond lemma" after the celebrated
Newman's Diamond Lemma for graphs \cite{newman}.

Shirshov's lemma is now called the Composition-Diamond lemma for Lie
and associative algebras.

 This kind of ideas were independently discovered by H. Hironaka
\cite{H64} for power series algebras and by B. Buchberger
\cite{bu65, bu70} for polynomial algebras. Buchberger suggested the
name ``Gr\"{o}bner bases". Applications of Gr\"{o}bner bases in
mathematics (particularly in algebraic geometry) and in computer
science are now well established.

At present, there are quite a few Composition-Diamond lemmas
(CD-lemma for short) for different classes of non-commutative and
non-associative algebras. Let us mention some.

A. I. Shirshov \cite{s62'} proved a CD-lemma for commutative
(anti-commutative) non-associative algebras, and mentioned that this
lemma is also valid for non-associative algebras. It was used to
solve the word problem for these classes of algebras. For
non-associative algebras, this (but not the CD-lemma) was known, see
A. I. Zhukov \cite{Z}.

 A.~A.~Mikhalev \cite{MikhPetr} proved a CD-lemma for Lie super-algebras.

 T. Stokes \cite{st90} proved a CD-lemma for left ideals of an algebra $\mathbbm{k}[X]\otimes E_\mathbbm{k}(Y)$,
 the tensor product of an exterier (Grassman) algebra and a polynomial algebra.

A. A. Mikhalev and E. A. Vasilieva \cite{MV} proved a CD-lemma for
the free supercommutative polynomial algebras.

 A. A. Mikhalev and A.~A. Zolotykh  \cite{MZ} proved a CD-lemma
 for $\mathbbm{k}[X]\otimes \mathbbm{k}\langle Y\rangle$, the tensor product of a polynomial algebra and a free algebra.

L. A. Bokut, Y. Fong and W. F. Ke \cite{bfk} proved a CD-lemma for
associative conformal algebras.

L. Hellstr\"{o}m \cite{H} proved a CD-lemma for a non-commutative
power series algebra.

 S.-J. Kang and K.-H. and Lee \cite{kl1, kl3} and E. S. Chibrikov \cite{ch} proved a
CD-lemma for a module over an algebra (see also \cite{cyz}).

D. R. Farkas, C. D. Feustel and E. L. Green \cite{ffg93} proved a
CD-lemma for path algebras.

L. A. Bokut and K. P. Shum \cite{b06} proved a CD-lemma for
$\Gamma$-algebras.

Y. Kobayashi \cite{kobayashi}  proved a CD-lemma for algebras based
on well-ordered semigroups,  L. A. Bokut, Yuqun Chen and Cihua Liu
\cite{bcl} proved a CD-lemma for dialgebras (see also \cite{bc}), L.
A. Bokut, Yuqun Chen and Yongshan Chen \cite{BCC} proved a CD-lemma
for tensor product of free algebras, and L. A. Bokut, Yuqun Chen and
Jianjun Qiu \cite{BCQ} proved a CD-lemma for associative algebras
with multiple linear operators.

Let $\mathscr{D}$ be a set of symbols. Let $\mathbbm{k}$ be a
commutative ring with unit and $\mathscr{A}$ an associative algebra
over $\mathbbm{k}$. Then a pair $(\mathscr{A},
\partial_{\mathscr{A}})$ is called a differential algebra with
differential operators $\mathscr{D}$ or $\mathscr{D}$-algebra for
short if $\partial_{\mathscr{A}}: \ \mathscr{D}\rightarrow
Der\mathscr{A}$ is an injective map, where $Der\mathscr{A}$ stands
for the set of all derivations of $\mathscr{A}$, i.e., any
$\delta\in Der\mathscr{A}, \ a,b \in \mathscr{A}$,
$$
\delta(ab)=\delta(a)\cdot b+a\cdot \delta(b).
$$
Given a $\mathscr{D}$-algebra $\mathscr{A}$ in this sense, we
identify $D\in \mathscr{D}$ with the operator
$\partial_{\mathscr{A}} D$.

In this paper, we establish a Composition-Diamond lemma for free
differential algebras. As applications, we give Gr\"{o}bner-Shirshov
bases for free Lie-differential algebras and free
commutative-differential algebras, respectively.

The authors would like to express their deepest gratitude to
Professor L. A. Bokut for his kind guidance, useful discussions and
enthusiastic encouragement. We also thank the referee for giving us
some useful remarks.

\section{Free differential algebras}

In this section, we construct a free differential algebra with
differential operators $\mathscr{D}$.

Let $\mathscr{A}$ be a $\mathscr{D}$-algebra. A subset $I$ of
$\mathscr{A}$ is a $\mathscr{D}$-ideal if $I$ is a
$\mathbbm{k}$-algebra ideal such that $D(I)\subseteq I$ for any
$D\in \mathscr{D}$.

Let $\mathscr{A},\mathscr{B}$ be $\mathscr{D}$-algebras and
$\psi:\mathscr{A}\rightarrow \mathscr{B}$ a map. Then $\psi$ is
called a $\mathscr{D}$-homomorphism if $\psi$ is a
$\mathbbm{k}$-algebra homomorphism such that for any $a\in
\mathscr{A} \ and \ D\in \mathscr{D}, \ \ \psi(D(a))=D(\psi(a))$,
which is precisely
$\psi(\partial_{\mathscr{A}}D(a))=\partial_{\mathscr{B}}D(\psi(a))$.

\begin{definition}
Let $X$ be a set, $\mathscr{D}(X)$ a $\mathscr{D}$-algebra and
$\iota:X\rightarrow\mathscr{D}(X)$ an inclusion map. If for any
$\mathscr{D}$-algebra $\mathscr{A}$ and any map $\varphi$:
$X\rightarrow\mathscr{A}$, there exists a unique
$\mathscr{D}$-homomorphism $\varphi^*:
\mathscr{D}(X)\rightarrow\mathscr{A}$ such that
$\varphi^*\iota=\varphi$, then $\mathscr{D}(X)$ is called a free
$\mathscr{D}$-algebra generated by $X$.
\end{definition}

Now we construct a free $\mathscr{D}$-algebra directly.

Let $\mathscr{D}=\{D_j|j\in J\}$ and $\mathbb{N}$ the set of
non-negative integers. For any
 $m\in \mathbb{N}$ and $\bar{j}=(j_1,\cdots,j_m)\in J^m$,
denote by $D^{\bar{j}}=D_{j_1}D_{j_2}\cdots D_{j_{m}}$ and
$D^{\omega}(X)=\{D^{\bar{j}}(x)|x\in X, \ \bar{j}\in J^m, \
m\in\mathbb{N}\}$, where $D^{\bar{j}}(x)=x$ if $\bar{j}\in J^0$. Let
$T=(D^{\omega}(X))^*$ be the free monoid generated by
$D^{\omega}(X)$. For any
$u=D^{\overline{i_1}}(x_1)D^{\overline{i_2}}(x_2)\cdots
D^{\overline{i_n}}(x_n)$ $\in T$, the length of $u$, denoted by
$|u|$, is defined to be $n$. In particular, $|1|=0$.

Let $\mathscr{D}(X)=\mathbbm{k}T$ be the $\mathbbm{k}$-algebra
spanned by $T$. For any $D_j\in \mathscr{D}$,
 we define the $\mathbbm{k}$-linear map $ D_j: \ \mathscr{D}(X) \rightarrow
\mathscr{D}(X)$ by induction on $|u|$ for $u\in T$:
\begin{enumerate}
\item[(1)] $D_j(1)=0$.
\item[(2)] Suppose that $u=D^{\bar{i}}(x)=D_{i_1}D_{i_2}\cdots
D_{i_{m_i}}(x)$. Then $D_j(u)=D_jD_{i_1}D_{i_2}\cdots
D_{i_{m_i}}(x)$.
\item[(3)] Suppose that $u=D^{\bar{i}}(x)\cdot v, \ v\in T$. Then
$D_j(u)=(D_jD^{\bar{i}}(x))\cdot v+D^{\bar{i}}(x)\cdot D_j(v)$.
\end{enumerate}
By definition, $D_j$ is a derivation on $\mathscr{D}(X)$ and then
$\mathscr{D}(X)$ is a differential algebra with differential
operators $\mathscr{D}$.

\begin{theorem}
$\mathscr{D}(X)$ is a free $\mathscr{D}$-algebra generated by $X$.
\end{theorem}

\noindent{\bf Proof.}
 Let $\mathscr{A}$ be any
$\mathscr{D}$-algebra and $\varphi:X\rightarrow\mathscr{A}$ a map.
We define a linear map
$\varphi^*:\mathscr{D}(X)\rightarrow\mathscr{A}$ by
$$
\varphi^*(D^{\overline{i_1}}(x_1)D^{\overline{i_2}}(x_2)\cdots
D^{\overline{i_k}}(x_k))=D^{\overline{i_1}}(\varphi(x_1))D^{\overline{i_2}}(\varphi(x_2))\cdots
D^{\overline{i_k}}(\varphi(x_k)).
$$
Then, it is easy to check that $\varphi^*$ is the unique
$\mathscr{D}$-homomorphism such that $\varphi^*\iota=\varphi$, where
$\iota:X\rightarrow\mathscr{D}(X)$ is the inclusion map and the
proof is complete.  $\ \ \ \square$

\section{Composition-Diamond lemma for  free differential algebras}

Let $\mathscr{D}=\{D_j|j\in J\}$, $X$ and $J$ well ordered sets,
 $D^{\bar{i}}(x)=D_{i_1}D_{i_2}\cdots D_{i_{m}}(x)\in D^{\omega}(X)$  and
$$
wt(D^{\bar{i}}(x))=(x; m, i_1, i_2,\cdots, i_{m}).
$$
We order $D^{\omega}(X)$ as follows:
$$
D^{\bar{i}}(x)> D^{\bar{j}}(y)\Longleftrightarrow
wt(D^{\bar{i}}(x))> wt(D^{\bar{j}}(y)) \ \mbox{ lexicographically}.
$$
It is easy to check that this order is a well ordering on
$D^{\omega}(X)$.

Now, we order $T=(D^{\omega}(X))^*$ by deg-lex order $>$, i.e.,  two
words are first compared by length and then lexicographically. We
will use this order in the sequel.

Then, for any $0\neq f\in \mathscr{D}(X)$, we can write
$$
f=\alpha_{\bar{f}}\bar{f}+\sum_{i=1}^{k}\alpha_iu_i,
$$
where $\bar{f},u_i\in T, \ 0\neq \alpha_{\bar{f}},\alpha_i\in
\mathbbm{k}$, and $\bar{f}>u_i$.  $\bar{f}$ is called the leading
term of $f$ and $\alpha_{\bar{f}}$ the leading coefficient. If
$\alpha_{\bar{f}}=1$, we say that $f$ is monic.

The proofs of Lemmas \ref{l1}, \ref{l2} and \ref{l3} are
straightforward.

\begin{lemma}\label{l1}
If $u>v$, $u,v\in T$, then for any $a,b\in T$, $aub>avb$. $ \square$
\end{lemma}

\begin{lemma}\label{l2}
Let $D^{\bar{i}}(x)=D_{i_1}D_{i_2}\cdots D_{i_{m}}(x)\in D^{w}(X), \
u=D^{\bar{i}}(x)v,\ v\in T$. Then, for any $D_j\in \mathscr{D}$, we
have
$$
\overline{D_j(u)}=D_jD_{i_1}D_{i_2}\cdots D_{i_{m}}(x)v. \ \ \ \
\square
$$
\end{lemma}

\begin{lemma}\label{l3}
For any $u,v\in T\backslash\{1\}$, if $u>v$, then
$\overline{D(u)}>\overline{D(v)}$ for any $D\in \mathscr{D}. \ \ \ \
\ \square$
\end{lemma}

By the above lemmas, the order $>$ on $T$ is monomial in the sense
that for any $u,v\in T$,
$$
u>v\Rightarrow aub>avb \ \mbox{ and } \
\overline{D(u)}>\overline{D(v)} \  \mbox{ for any } a,b\in T,\
D\in\mathscr{D},
$$
where we assume that $v\neq 1$ for the latter case.

\begin{corollary}\label{l4}
Let $D^{\bar{i}}(x)=D_{i_1}D_{i_2}\cdots D_{i_{m}}(x)\in D^{w}(X), \
u=D^{\bar{i}}(x)v,\ v\in T$. Then, for
$D^{\bar{j}}=D_{j_1}D_{j_2}\cdots D_{j_{n}}$,
$$
\overline{D^{\bar{j}}(u)}=D_{j_1}D_{j_2}\cdots
D_{j_{n}}D_{i_1}D_{i_2}\cdots D_{i_{m}}(x)v. \ \ \square
$$
\end{corollary}

Let the notation be as in Corollary \ref{l4}. For convenience, we
denote $\overline{D^{\bar{j}}(u)}$ by ${d^{\bar{j}}(u)}$.

\begin{definition}
Let $f,g\in \mathscr{D}(X)$ be monic polynomials and $w,a,b\in T$.
Then there are two kinds of compositions.
\begin{enumerate}
\item[(i)] There are two sorts of composition of inclusion:

If $w=\bar{f}=a\cdot d^{\bar{j}}(\bar{g})\cdot b$, then the
composition is
$$
(f,g)_{w}=f-a\cdot D^{\bar{j}}(g)\cdot b.
$$
If $w=d^{\bar{i}}(\bar{f})=\bar{g}\cdot b$, then the composition is
$$
(f,g)_{w}=D^{\bar{i}}(f)-g\cdot b.
$$

\item[(ii)] Composition of intersection:

If $w=\bar{f}\cdot b=a\cdot d^{\bar{j}}(\bar{g})$ such that
$|\bar{f}|+|\bar{g}|>|w|$, then the composition is
$$
(f,g)_{w}=f\cdot b-a\cdot D^{\bar{j}}(g).
$$
In this case, we assume that $a,b\neq 1$.
\end{enumerate}
\end{definition}

By Lemmas \ref{l1} and \ref{l3}, we have

\begin{lemma}\label{l0}
For any composition $(f,g)_{w}$, we have $\overline{(f,g)_{w}}<w.$ \
\ $\square$
\end{lemma}

Let $S\subset \mathscr{D}(X)$ be a subset of monic polynomials. Then
we define $(S,\mathscr{D})$-words $u_s$ inductively: for any $s\in
S$,
\begin{enumerate}
\item[(i)] for any $D^{\bar{j}}=D_{j_1}D_{j_2}\cdots D_{j_{m}}, \ u_s=D^{\bar{j}}(s)$ is
an $(S,\mathscr{D})$-word of $(S,\mathscr{D})$-length 1;
\item[(ii)] if $u_s$ is an $(S,\mathscr{D})$-word of $(S,\mathscr{D})$-length $k$ and $v$
is a word in $T$ of length $l$, then
\begin{equation*}
u_s\cdot v \ \  and \ \  v\cdot u_s
\end{equation*}
are $(S,\mathscr{D})$-words of $(S,\mathscr{D})$-length $k+l$.
\end{enumerate}

Then, any element of $Id(S)$, the $\mathscr{D}$-ideal of
$\mathscr{D}(X)$ generated by $S$, is a linear combination of
$(S,\mathscr{D})$-words.

By Lemmas \ref{l1}, \ref{l3} and \ref{l4}, we obtain

\begin{lemma} For any $(S,\mathscr{D})$-word $u_s=a\cdot D^{\bar{j}}(s)\cdot b, \ a,b\in T,\ s\in
S$, we have  $\overline{u_s}=a\cdot d^{\bar{j}}(\bar{s})\cdot b$. \
\  $\square$ \end{lemma}

\begin{definition}Let $S\subset\mathscr{D}(X) $ be a monic subset and $w\in T$.
Then a polynomial $f\in\mathscr{D}(X)$ is said to be trivial modulo
$(S,w)$, denoted by $f\equiv0 \ mod(S,w)$, if
$$
f=\sum_i\alpha_iu_{s_i},
$$
where each $u_{s_i}$ is an $(S,\mathscr{D})$-word, $\alpha_i\in
\mathbbm{k}$ and $\overline{u_{s_i}}<w$.

$S$ is called a Gr\"{o}bner-Shirshov basis in $\mathscr{D}(X)$ if
all compositions of elements of $S$ are trivial modulo $S$ and
corresponding $w$.
\end{definition}

\begin{lemma}\label{la}
$S$ is a Gr\"{o}bner-Shirshov basis in $\mathscr{D}(X)$ if and only
if for any $f,g\in S$,
\begin{enumerate}
\item[(i)] if $w=d^{\bar{i}}(\bar{f})=a\cdot d^{\bar{j}}(\bar{g})\cdot b$ for
some $w,a,b\in T$, then
$$
D^{\bar{i}}(f)-a\cdot D^{\bar{j}}(g)\cdot b\equiv0 \ \ \ \ mod(S,w);
$$

\item[(ii)] if $w=d^{\bar{i}}(\bar{f})\cdot b=a\cdot d^{\bar{j}}(\bar{g})$
for some $w,a,b\in T, a,b\neq1$ such that $|\bar{f}|+|\bar{g}|>|w|$,
then
$$
D^{\bar{i}}(f)\cdot b-a\cdot D^{\bar{j}}(g)\equiv0 \ \ \ \ mod(S,w).
$$
\end{enumerate}
\end{lemma}
\noindent{\bf Proof.} \ The sufficiency is clear. We only need to
prove the necessity.

Firstly,  $w=d^{\bar{i}}(\bar{f})=a\cdot d^{\bar{j}}(\bar{g})\cdot
b$ for some $w,a,b\in T$.

\emph{Case 1.1} If $a\neq1$, then $\exists  \ a'\in T$, such that
$d^{\bar{i}}a'=a$. We have $w'=\bar{f}=a'd^{\bar{j}}(\bar{g})\cdot
b$, $(f,g)_{w'}=f-a'\cdot D^{\bar{j}}(g)\cdot b\equiv0 \ \ \
mod(S,w')$ and $w=d^{\bar{i}}w'$. Therefore,
\begin{eqnarray*}
&&D^{\bar{i}}(f)-a\cdot D^{\bar{j}}(g)\cdot b\\
&=&D^{\bar{i}}(f-a'\cdot D^{\bar{j}}(g)\cdot b)+D^{\bar{i}}(a'\cdot
D^{\bar{j}}(g)\cdot
b)-a\cdot D^{\bar{j}}(g)\cdot b \\
&=&D^{\bar{i}}(f-a'\cdot D^{\bar{j}}(g)\cdot b)+D^{\bar{i}}(a')\cdot
D^{\bar{j}}(g)\cdot b+a'\cdot D^{\bar{i}}(D^{\bar{j}}(g)\cdot
b)-a\cdot D^{\bar{j}}(g)\cdot b \\
&=&D^{\bar{i}}(f-a'\cdot D^{\bar{j}}(g)\cdot b)+d^{\bar{i}}(a')\cdot
D^{\bar{j}}(g)\cdot b+ \sum_ic_i\cdot D^{\bar{j}}(g)\cdot b +a'\cdot
D^{\bar{i}}(D^{\bar{j}}(g)\cdot
b)-a\cdot D^{\bar{j}}(g)\cdot b \\
&=&D^{\bar{i}}((f,g)_{w'})+ \sum_kc_k\cdot D^{\bar{j}}(g)\cdot b
+a'\cdot D^{\bar{i}}(D^{\bar{j}}(g)\cdot
b) \\
\end{eqnarray*}
where $D^{\bar{i}}(a')=d^{\bar{i}}(a)+\sum_kc_k=a+\sum_kc_k, \
c_k<a$ and all $c_k\cdot D^{\bar{j}}(g)\cdot b, \ a'\cdot
D^{\bar{i}}(D^{\bar{j}}(g)\cdot b)$ are $(S,\mathscr{D})$-words with
the leading terms less than $w$. So we get
$$
D^{\bar{i}}(f)-a\cdot D^{\bar{j}}(g)\cdot b\equiv0 \ \ \ mod(S,w).
$$

The following three cases can be proved similarly to the Case 1.1
and we omit the details.

\emph{Case 1.2} $a=1$ and $\exists \ \bar{l}$ such that
$\bar{f}=d^{\bar{l}} ( \bar{g})b$, $d^{\bar{j}} (
\bar{g})=d^{\bar{i}}(d^{\bar{l}}(\bar{g}))$. Then
$w'=\bar{f}=d^{\bar{l}}(\bar{g})\cdot b$, $w=d^{\bar{i}}w'$ and
$(f,g)_{w'}=f-D^{\bar{l}}(g)\cdot b\equiv0 \ \ \ mod(S,w')$.

\emph{Case 1.3} $a=1$ and $\exists \ \bar{l}$ such that $d^{\bar{l}}
(\bar{f})= \bar{g}b$, $d^{\bar{j}}(d^{\bar{l}} (
\bar{f}))=d^{\bar{i}}(\bar{f})$. Then
$w'=d^{\bar{l}}(\bar{f})=\bar{g}\cdot b$, $w=d^{\bar{j}}w'$ and
$(f,g)_{w'}=D^{\bar{l}} (f)-gb\equiv0 \ \ \ mod(S,w')$.

Secondly, $w=d^{\bar{i}}(\bar{f})\cdot b=a\cdot
d^{\bar{j}}(\bar{g})$ for some $w,a,b\in T, a,b\neq1$ such that
$|\bar{f}|+|\bar{g}|>|w|$. Then there exists an $a'\in T, \ a'\neq1$
such that $a=d^{\bar{i}}(a')$, $w'=\bar{f}\cdot b=a'\cdot
d^{\bar{j}}(\bar{g})$, $w=d^{\bar{i}}w'$ and $(f,g)_{w'}=f\cdot
b-a'\cdot D^{\bar{j}}(g)\equiv0 \ \ \ mod(S,w')$. $\square$

\begin{lemma}\label{l5}
Let $S$ be a Gr\"{o}bner-Shirshov basis in $\mathscr{D}(X)$ and
$s_1,s_2\in S$. If
$$
w=a_1\cdot d^{\bar{i}}(\bar{s_1})\cdot b_1=a_2\cdot
d^{\bar{j}}(\bar{s_2})\cdot b_2
$$
for some
$a_1,a_2,b_1,b_2\in T$, then
$$
f=a_1\cdot D^{\bar{i}}(s_1)\cdot b_1-a_2\cdot D^{\bar{j}}(s_2)\cdot
b_2\equiv0 \ \ mod(S,w).
$$
\end{lemma}

\noindent{\bf Proof.} \ There are three cases to consider.

\emph{Case 1}. $d^{\bar{i}}(\bar{s}_1)$ and $d^{\bar{j}}(\bar{s}_2)$
are mutually disjoint. We may assume that $d^{\bar{i}}(\bar{s}_1)$
is at the left of $d^{\bar{j}}(\bar{s}_2)$, i.e.,
$$
a_2=a_1\cdot d^{\bar{i}}(\bar{s}_1)\cdot a, \ \ \ b_1=a \cdot
d^{\bar{j}}(\bar{s}_2)\cdot b_2
$$
where $a\in T$. Then
\begin{eqnarray*}
f&=&a_1\cdot D^{\bar{i}}(s_1)\cdot b_1-a_2\cdot
D^{\bar{j}}(s_2)\cdot b_2\\
&=&a_1\cdot D^{\bar{i}}(s_1)\cdot a\cdot d^{\bar{j}}(\bar{s}_2)\cdot
b_2-a_1\cdot
d^{\bar{i}}(\bar{s}_1)\cdot a\cdot D^{\bar{j}}(s_2)\cdot b_2\\
&=&a_1\cdot(D^{\bar{i}}(s_1)-d^{\bar{i}}(\bar{s}_1))\cdot a\cdot
D^{\bar{j}}(s_2)\cdot b_2-a_1\cdot D^{\bar{i}}(s_1)\cdot a(
D^{\bar{j}}(s_2)-d^{\bar{j}}(\bar{s}_2))\cdot b_2\\
&=&\sum_l \alpha_la_l'\cdot D^{\bar{j}}(s_2)\cdot b_2-\sum_k
\beta_ka_1\cdot D^{\bar{i}}(s_1)\cdot b_k'
\end{eqnarray*}
where $\alpha_l,\beta_k\in \mathbbm{k}, \ a_l',b_k'\in T, \
a_l'<a_1\cdot d^{\bar{i}}(\bar{s}_1)\cdot a=a_2, \ b_k'<a \cdot
d^{\bar{j}}(\bar{s}_2)\cdot b_2=b_1$. Thus, $f\equiv0 \ \ mod(S,w).$

\emph{Case 2}. One of $d^{\bar{i}}(\bar{s}_1)$ and
$d^{\bar{j}}(\bar{s}_2)$ is a subword of the other, say,
$w'=d^{\bar{i}}(\bar{s}_1)=a\cdot d^{\bar{j}}(\bar{s}_2)\cdot b$.
Then $a_2=a_1a, \ b_2=bb_1$ and
\begin{eqnarray*}
f&=&a_1\cdot D^{\bar{i}}(s_1)\cdot b_1-a_2\cdot
D^{\bar{j}}(s_2)\cdot b_2\\
&=&a_1\cdot D^{\bar{i}}(s_1)\cdot b_1- a_1a\cdot
D^{\bar{j}}(s_2)\cdot bb_1\\
&=&a_1\cdot( D^{\bar{i}}(s_1)-a\cdot D^{\bar{j}}(s_2)\cdot
b)\cdot b_1\\
&\equiv&0 \  \ \ \ \ \ \ \ \  \ \ \ \ \ \ \ \ \ \ \ \
mod(S,w)
\end{eqnarray*}
since  $w=a_1w' b_1$ and $ \ D^{\bar{i}}(s_1)-a\cdot
D^{\bar{j}}(s_2)\cdot b\equiv0 \ \ mod(S,w')$ by Lemma \ref{la}.

 \emph{Case 3}. $d^{\bar{i}}(\bar{s}_1)$
and $d^{\bar{j}}(\bar{s}_2)$ have a nonempty intersection as a
subword of $w$, but $d^{\bar{i}}(\bar{s}_1)$ and
$d^{\bar{j}}(\bar{s}_2)$ are not subwords of each other. We may
assume that
$$
a_2=a_1a, \ b_1=bb_2, w'=d^{\bar{i}}(\bar{s}_1)\cdot b=a\cdot
d^{\bar{j}}(\bar{s}_2).
$$
Then, we have
\begin{eqnarray*}
f&=&a_1\cdot D^{\bar{i}}(s_1)\cdot b_1-a_2\cdot
D^{\bar{j}}(s_2)\cdot b_2\\
&=&a_1\cdot D^{\bar{i}}(s_1)\cdot bb_2- a_1a\cdot
D^{\bar{j}}(s_2)\cdot b_2\\
&=&a_1\cdot( D^{\bar{i}}(s_1)\cdot b-a\cdot
D^{\bar{j}}(s_2))\cdot b_2\\
&\equiv&0 \  \ \ \ \ \ \ \ \  \ \ \ \ \ \ \ \ \ \ \
mod(S,w)
\end{eqnarray*}
since $w=a_1w'b_2$ and also by Lemma \ref{la},
$D^{\bar{i}}(s_1)\cdot b-a\cdot D^{\bar{j}}(s_2)\equiv 0 \ \
mod(S,w')$.

This completes the proof. $ \square$

\begin{lemma}\label{l6}
Let $S\subset \mathscr{D}(X)$ be a monic subset and $Irr(S)=\{ u\in
T \ | \ u\neq a\cdot d^{\bar{i}}(\bar{s})\cdot b  \ for \ all \ s\in
S, \  a,b\in T, \ \bar{i}\in J^m, \ m\in \mathbb{N}\}$. Then for any
$f\in \mathscr{D}(X)$,
\begin{equation*}
f=\sum\limits_{u_i\leq \bar f }\alpha_iu_i+
\sum\limits_{\overline{v_{s_j}}\leq \bar f}\beta_jv_{s_j}
\end{equation*}
where each $\alpha_i,\beta_j\in \mathbbm{k}, \ u_i \in Irr(S)$ and
$v_{s_j}$ is an $(S,\mathscr{D})$-word.
\end{lemma}
{\bf Proof.} Let
$f=\sum\limits_{i}\alpha_{i}u_{i}\in{\mathscr{D}(X)}$, where
$0\neq{\alpha_{i}\in{\mathbbm{k}}}$ and $u_{1}>u_{2}>\cdots$. If
$u_1\in{Irr(S)}$, then let $f_{1}=f-\alpha_{1}u_1$. If
$u_1\not\in{Irr(S)}$, then there exist some $s_1\in{S},\ \bar{i_1}$
and $a_1,b_1\in{T}$, such that $\bar
f=u_1=a_1d^{\bar{i_1}}(\bar{s_1})b_1$. Let
$f_1=f-\alpha_1a_1D^{\bar{i_1}}(s_1)b_1$. In both cases, we have
$\bar{f_1}<\bar{f}$. Then the result follows by induction on
$\bar{f}$. \ \ \ \ $\square$

\ \

The following theorem is an analogue  of the Shirshov's Composition
lemma for Lie algebras \cite{s62}, which was specialized to
associative algebras by Bokut \cite{b76}, see also Bergman \cite{b}.

\begin{theorem}\label{3}
\textbf{(Composition-Diamond lemma for differential algebras)} Let
$\mathscr{D}(X)$ be the free differential algebra with differential
operators $\mathscr{D}=\{D_j|
 j\in J\}$,  $S\subset \mathscr{D}(X)$ a
monic subset, $Id(S)$ the $\mathscr{D}$-ideal of $\mathscr{D}(X)$
generated by $S$ and  $<$ the order on $T=(D^{\omega}(X))^*$ as
before. Then the following statements are equivalent:
\begin{enumerate}
\item[(i)]$S$ is a  Gr\"{o}bner-Shirshov basis in $\mathscr{D}(X)$.

\item[(ii)] $f\in Id(S)\Rightarrow\bar{f} =a\cdot
d^{\bar{i}}(\bar{s})\cdot b $ for some $s\in S,\ \bar{i}\in J^m, \
m\in \mathbb{N}$ and $a,b\in T.$

\item[(iii)] $Irr(S)=\{ u\in T \ | \ u\neq a\cdot
d^{\bar{i}}(\bar{s})\cdot b  \ for \ all \ s\in S, \  a,b\in T, \
\bar{i}\in J^m, \ m\in \mathbb{N}\}$ is a $\mathbbm{k}$-linear basis
of $\mathscr{D}(X|S)=\mathscr{D}(X)/Id(S)$.
\end{enumerate}
\end{theorem}

\noindent{\bf Proof.} \ \ $(i)\Rightarrow (ii)$. Let $S$ be a
Gr\"{o}bner-Shirshov basis and $0\neq f\in Id(S).$ Then

$$f=\sum_{i=1}^{k}\alpha_ia_i\cdot D^{\overline{j_{i}}}(s_i)\cdot b_i, $$
\noindent where $\alpha_i\in \mathbbm{k}, \ a_i,b_i\in T, \ s_i\in
S$. Let $w_i=a_i\cdot d^{\overline{j_{i}}}(\bar{s_i})\cdot b_i.$  We
may assume without loss of generality that
$$w_1=w_2=\cdots =w_l>w_{l+1}\geq w_{l+2}\geq \cdots \geq w_k$$
\noindent for some $l\geq 1$. We proceed by induction on $w_1$ and
$l$.

If $l=1$, then $\bar{f}=w_1=a_1\cdot d^{\bar{j_1}}(\bar{s_1})\cdot
b_1$ and we are done.

If $\alpha_1+\alpha_2\neq 0$ or $l>2$, we have $a_1\cdot
d^{\bar{j_1}}(\bar{s_1})\cdot b_1=w_1=w_2=a_2\cdot
d^{\bar{j_2}}(\bar{s_2})\cdot b_2$ and furthermore,
\begin{eqnarray*}
&&\alpha_1a_1\cdot D^{\bar{j_1}}(s_1)\cdot b_1+\alpha_2a_2\cdot D^{\bar{j_2}}(s_2)\cdot b_2\\
&=&(\alpha_1+\alpha_2)a_1\cdot D^{\bar{j_1}}(s_1)\cdot
b_1+\alpha_2(a_2\cdot D^{\bar{j_2}}(s_2)\cdot b_2-a_1\cdot
D^{\bar{j_1}}(s_1)\cdot b_1).
\end{eqnarray*}

\noindent By Lemma \ref{l5}, the last item of the above is a linear
combination of  $(S,\mathscr{D})$-words less than $w_1$. The result
follows by induction on $l$.

For the case of $\alpha_1+\alpha_2= 0$ and $l=2$, we use induction
on $w_1$ and the result follows.
\\

$(ii)\Rightarrow (iii).$ For any $f\in \mathscr{D}(X)$, by Lemma
\ref{l6}, we have
\begin{equation*}
f=\sum\limits_{u_i\leq \bar f }\alpha_iu_i+
\sum\limits_{\overline{v_{s_j}}\leq \bar f}\beta_jv_{s_j},
\end{equation*}
where each $\alpha_i,\beta_j\in \mathbbm{k}, \ u_i \in Irr(S)$ and
$v_{s_j}$ is an $(S,\mathscr{D})$-word. Therefore
$$
f+Id(S)=\sum_i\alpha_i(u_i+Id(S)).
$$ On the other hand,
suppose that $\sum\limits_{i}\alpha_iu_i=0$ in
$\mathscr{D}(X)/Id(S)$, where $\alpha_i\in \mathbbm{k}$, $u_i\in
{Irr(S)}$. This means that $\sum\limits_{i}\alpha_iu_i\in{Id(S)}$ in
$\mathscr{D}(X)$. Then all $\alpha_i$ must be equal to zero.
Otherwise, $\overline{\sum\limits_{i}\alpha_iu_i}=u_j\in{Irr(S)}$
for some $j$ which contradicts (ii).
 \\

 $(iii)\Rightarrow (i).$
 For any $f,g\in{S}$ , by Lemmas \ref{l0} and  \ref{l6}, we
have
$$
(f,g)_{w}=\sum\limits_{u_i\in Irr(S),\ u_i<w }\alpha_iu_i+
\sum\limits_{\overline{v_{s_j}}< w}\beta_jv_{s_j}.
$$
Since $(f,g)_{w}\in {Id(S)}$ and by $(iii)$, we have
$$
(f,g)_{w}= \sum\limits_{\overline{v_{s_j}}< w}\beta_jv_{s_j}.
$$
Therefore, $S$ is a Gr\"{o}bner-Shirshov basis. \ \ \ \ $\square$

\section{Gr\"{o}bner-Shirshov bases for free Lie-differential algebras}

In this section, we define Lie- (commutative-) differential
algebras. As applications of the Composition-Diamond lemma for
differential algebras (Theorem \ref{3}), we give
Gr\"{o}bner-Shirshov bases for free Lie-differential algebras and
free commutative-differential algebras, respectively.

It is well-known that for an arbitrary algebra $\mathscr{A}$ the set
$Der\mathscr{A}$ is a Lie algebra with respect to the ordinary
commutator of linear maps.
\begin{definition}
Let $(\mathscr{D},[,])$ be a Lie algebra over $\mathbbm{k}$. A
Lie-differential algebra (Lie-$\mathscr{D}$ algebra for short) is a
$\mathscr{D}$-algebra $(\mathscr{A},\partial)$, where
$\partial:\mathscr{D}\rightarrow Der\mathscr{A}$ is a homomorphism
of Lie algebras.
\end{definition}

\begin{definition}
Let $X$ be a set, $(\mathscr{D}, [,])$ a Lie algebra,
$L\mathscr{D}(X)$ a Lie-$\mathscr{D}$ algebra and
$\iota:X\rightarrow L\mathscr{D}(X)$ an inclusion map. If for any
Lie-$\mathscr{D}$ algebra $\mathscr{A}$ and any map $\varphi$:
$X\rightarrow\mathscr{A}$, there exists a unique
$\mathscr{D}$-homomorphism $\varphi^*: \
L\mathscr{D}(X)\rightarrow\mathscr{A}$ such that
$\varphi^*\iota=\varphi$, then $L\mathscr{D}(X)$ is called a free
Lie-$\mathscr{D}$ algebra generated by $X$.
\end{definition}

The proof of the following theorem is straightforward and we omit
the details.

\begin{theorem}\label{t4.3}
Let $(\mathscr{D}, [,])$ be a Lie algebra with a $\mathbbm{k}$-basis
$\{D_j \ | \ j\in J\}$ and $[D_i,D_j]=D_iD_j-D_jD_i=\sum_r
\alpha_{ij}^rD_r$, $i,j\in J$, where each $\alpha_{ij}^r\in
\mathbbm{k}$. Let $\mathscr{D}(X)$ be a free $\mathscr{D}$-algebra
and
\begin{equation*}
S=\{D_iD_j(u)-D_jD_i(u)-\sum_r \alpha_{ij}^rD_r(u) \ | \ u\in T, \
i,j\in J\}.
\end{equation*}
Then  $\mathscr{D}(X|S)=\mathscr{D}(X)/Id(S)$ is a free
Lie-$\mathscr{D}$ algebra generated by $X$. \ \ \ \ $\square$
\end{theorem}

Let
\begin{eqnarray*}
&&H=\{D_{i_1}\cdots D_{i_{m}}(x)\ |\  i_1\leq\cdots \leq i_{m}, \
i_1, \cdots, i_m\in J, \
m\in\mathbb{ N}, \ x\in X\},\\
&&D^{\bar{i}}(x)=D_{i_1}\cdots D_{i_{m}}(x),\\
 &&S_0=\{D_pD_qD^{\bar{i}}(x)-D_qD_pD^{\bar{i}}(x)-
 \sum_r\alpha_{pq}^rD_rD^{\bar{i}}(x) \ |
  \  p>q, \  p,q\in J, \ D^{\bar{i}}(x)\in H, \ x\in X\}.
\end{eqnarray*}

\begin{lemma}\label{l4.4}
$Irr(S_0)=H^*$, where $H^*$ is the free monoid generated by $H$. \ \
$\square$
\end{lemma}

\begin{lemma}\label{l3.6}
Suppose $S$ and $S_0$ are sets defined as above. Then, they generate
the same $\mathscr{D}$-ideal in $\mathscr{D}(X)$, i.e.,
\begin{equation*}
Id(S)=Id(S_0).
\end{equation*}

\end{lemma}
{\bf Proof.} Since $S_0$ is a subset of $S$, it suffices to prove
that $\mathscr{D}(X|S_0)$ is a Lie-$\mathscr{D}$  algebra. We need
only prove that in $\mathscr{D}(X|S_0)$, for any $u\in T$,
\begin{equation*}
D_iD_j(u)-D_jD_i(u)-\sum_r \alpha_{ij}^rD_r(u)=0.
\end{equation*}
By Lemma \ref{l6}, it suffices to prove that for any $v\in
Irr(S_0)$,
\begin{equation*}
D_iD_j(v)-D_jD_i(v)-\sum_r \alpha_{ij}^rD_r(v)=0.
\end{equation*}
We use induction on the length of $v$.

If $|v|=1$, then the result is trivial.

If $|v|\geq2$, i.e., $v=D^{\bar{t}}(x)v_1$ where $D^{\bar{t}}(x)\in
H, \ 1\neq v_1\in H^*$, then by induction,
\begin{eqnarray*}
& \ &D_iD_j(v)-D_jD_i(v)-\sum_r \alpha_{ij}^rD_r(v)
\\&=&D_iD_j(D^{\bar{t}}(x)v_1)-D_jD_i(D^{\bar{t}}(x)v_1)-\sum_r
\alpha_{ij}^rD_r(D^{\bar{t}}(x)v_1)\\
&=&(D_iD_jD^{\bar{t}}(x)-D_jD_iD^{\bar{t}}(x)-\sum_r
\alpha_{ij}^rD_rD^{\bar{t}}(x))v_1\\
& \ &+D^{\bar{t}}(x)(D_iD_j(v_1)-D_jD_i(v_1)-\sum_r
\alpha_{ij}^rD_r(v_1))\\
&=&0 \ . \ \ \ \ \square
\end{eqnarray*}

By Lemma \ref{l6}, we have

\begin{lemma}\label{l3.7}
For any $p,q\in J$ and $D^{\bar{i}}x\in D^{\omega}(X)$, if $p>q$,
then
$$
D_pD_qD^{\bar{i}}(x)-D_qD_pD^{\bar{i}}(x)-\sum_r\alpha_{pq}^rD_rD^{\bar{i}}(x)\equiv0
\ \ \ \ \ \ mod(S_0,w)
$$
where $w\in D^{\omega}(X)$ and $w>D_pD_qD^{\bar{i}}x$.\ \ $\square$
\end{lemma}

\begin{theorem}
Let $(\mathscr{D}, [,])$ be a Lie algebra with a $\mathbbm{k}$-basis
$\{D_j \ | \ j\in J\}$ and $[D_i,D_j]=D_iD_j-D_jD_i=\sum_r
\alpha_{ij}^rD_r$, $i,j\in J$, where each $\alpha_{ij}^r\in
\mathbbm{k}$. Let $\mathscr{D}(X)$ be the free differential algebra
with differential operators $\mathscr{D}$. Let the order $<$ on $T$
be as before and
\begin{eqnarray*}
 S_0&=&\{D_pD_qD_{i_1}\cdots D_{i_m}(x)-D_qD_pD_{i_1}\cdots D_{i_m}(x)-
\sum_r \alpha_{pq}^rD_rD_{i_1}\cdots D_{i_m}(x) \ |
 \\
&& \ \ \ p>q,\ i_1\leq \cdots\leq i_m, \ p,q,i_1, \cdots, i_m\in J,
\ m\geq0,  \ x\in X\}.
\end{eqnarray*}
Then
\begin{enumerate}
\item[(i)] $S_0$ is a Gr\"obner-Shirshov basis in
$\mathscr{D}(X)$.
\item[(ii)] $H^*$ is a $\mathbbm{k}$-basis of the free Lie-$\mathscr{D}$  algebra
$\mathscr{D}(X|S)$ generated by $X$, where $H=\{D_{i_1}\cdots
D_{i_{m}}(x)\ |\  i_1\leq\cdots \leq i_{m}, \ i_1, \cdots, i_m\in J,
\ m\in \mathbb{N}, \ x\in X\}$ and $H^*$ is the free monoid
generated by $H$.
\end{enumerate}
\end{theorem}
{\bf Proof.} (i). The possible composition is only the case of
inclusion $(f,g)_w=f-D_p g$ where
\begin{eqnarray*}
&&f=D_pD_qD_{i_1}\cdots D_{i_m}(x)-D_qD_pD_{i_1}\cdots D_{i_m}(x)-
\sum \alpha_{pq}^rD_rD_{i_1}\cdots D_{i_m}(x),\\
&&g=D_qD_{i_1}D_{i_2}\cdots D_{i_m}(x)-D_{i_1}D_qD_{i_2}\cdots
D_{i_m}(x)- \sum \alpha_{q{i_1}}^rD_rD_{i_2}\cdots D_{i_m}(x),\\
&&w=\bar{f}=\overline{D_p g}, \ p>q>i_1, \ i_1\leq i_2\leq
\cdots\leq i_m,\  m\geq 1.
\end{eqnarray*}

Now, we prove that $(f,g)_w\equiv 0 \ \ mod(S_0,w)$. Denote by
$u=D_{i_2}\cdots D_{i_m}(x)$. Then
\begin{align*}
(f,g)_w=D_pD_{i_1}D_qu+ \sum_r \alpha_{q{i_1}}^rD_pD_ru
-D_qD_pD_{i_1}u - \sum_r \alpha_{pq}^rD_rD_{i_1}u.
\end{align*}

Let $ A=D_pD_{i_1}D_qu, \  B=-D_qD_pD_{i_1}u $ and $ C=\sum_r
\alpha_{qi_1}^rD_pD_ru- \sum_r \alpha_{pq}^rD_rD_{i_1}u$. Then, by
Lemma \ref{l3.7}, we have
\begin{align*}
A&=D_pD_{i_1}D_qu \\
&\equiv D_{i_1}D_pD_qu+\sum_r \alpha_{p{i_1}}^rD_rD_{q}u\\
&\equiv D_{i_1}D_qD_pu+\sum_r \alpha_{pq}^rD_{i_1}D_ru+\sum_r
\alpha_{p{i_1}}^rD_rD_{q}u \ \ \ \ \  mod(S_0,w)
\end{align*}
and
\begin{align*}
B&=-D_qD_pD_{i_1}u \\
&\equiv -D_qD_{i_1}D_pu-\sum_r \alpha_{p{i_1}}^rD_{q}D_ru\\
&\equiv -D_{i_1}D_qD_pu-\sum_r \alpha_{q{i_1}}^rD_rD_pu-\sum_r
\alpha_{p{i_1}}^rD_qD_{r}u \ \ \ \ \  mod(S_0,w).
\end{align*}
Since $(\mathscr{D}, [,])$ is a Lie algebra, we have
$$
[D_i,D_j]=-[D_j,D_i] \ \  and \ \
[[D_i,D_j],D_k]+[[D_k,D_i],D_j]+[[D_j,D_k],D_i]=0.
$$
It follows that
$$\alpha_{ij}^s=-\alpha_{ji}^s \ \ and \ \
\sum_{s,t}  \alpha_{ij}^s\alpha_{sk}^t+\sum_{s,t}
\alpha_{ki}^s\alpha_{s{j}}^t+\sum_{s,t}
\alpha_{j{k}}^s\alpha_{si}^t=0.
$$
Thus,
\begin{align*}
(f,g)_w&= A+B+C \\
&\equiv\sum_r
\alpha_{pq}^rD_{i_1}D_ru+\sum_r\alpha_{p{i_1}}^rD_rD_{q}u
+\sum_r \alpha_{q{i_1}}^rD_pD_ru \\
& \ \ \ -\sum_r \alpha_{q{i_1}}^rD_rD_pu
-\sum_r\alpha_{p{i_1}}^rD_qD_{r}u-\sum_r
\alpha_{pq}^rD_rD_{i_1}u\\
& \equiv \sum_{r,l} \alpha_{pq}^r\alpha_{{i_1}r}^lD_lu+\sum_{r,l}
\alpha_{p{i_1}}^r\alpha_{r{q}}^lD_lu+\sum_{r,l}
\alpha_{q{i_1}}^r\alpha_{pr}^lD_lu\\
&\equiv(\sum_{r,l} \alpha_{pq}^r\alpha_{{i_1}r}^l+\sum_{r,l}
\alpha_{p{i_1}}^r\alpha_{r{q}}^l+\sum_{r,l}
\alpha_{q{i_1}}^r\alpha_{pr}^l)D_lu\\
 &\equiv0 \ \ \ \ \  mod(S_0,w).
\end{align*}
This shows (i).

(ii) follows from Theorem \ref{3}, Theorem \ref{t4.3}, Lemma
\ref{l4.4} and Lemma \ref{l3.6}.  \ \ \ \ $\square$

\ \

\noindent{\bf Remark}: Let $\mathscr{A}$ be a differential algebra
with differential operators $\mathscr{D}$. If $\mathscr{D}$ is a
commutative $\mathbbm{k}$-algebra, then $\mathscr{A}$ is called a
commutative-differential algebra. Clearly, a
commutative-differential algebra is a special Lie-differential
algebra. Therefore, all results in this section are true for
commutative-differential algebras.

\ \


\begin{thebibliography}{4}

\bibitem{b}G. M. Bergman,``The Diamond Lemma for Ring Theory", {\it Adv. in
Math.}, \textbf{29}(1978), p.178.


\bibitem{b76}L. A. Bokut, ``Imbeddings into Simple Associative
Algebras", {\it Algebra i Logika}, \textbf{15}(1976), p.117.



\bibitem{bc} L. A. Bokut, Yuqun Chen, ``Gr\"{o}bner-Shirshov
 Bases: some new Results", {\it Proceedings of the Second International Congress
 in Algebra and Combinatorics}, World Scientific, 2008, p.35.

\bibitem{BCC} L. A. Bokut, Yuqun Chen and Yongshan Chen, ``Composition-Diamond Lemma
for Tensor Product of Free Algebras", {\it J. Algebra}, to appear.
arXiv.org/abs/0804.2115.


\bibitem{bcl} L. A. Bokut, Yuqun Chen and Cihua Liu, ``Gr\"{o}bner-Shirshov
 Bases for Dialgebras", arXiv.org/abs/0804.0638.

\bibitem{BCQ} L. A. Bokut, Yuqun Chen and Jianjun Qiu, ``Gr\"{o}bner-Shirshov Bases for Associative Algebras with
Multiple Operators and Free Rota-Baxter Algebras",
arXiv.org/abs/0805.0640.

\bibitem{bfk}L. A. Bokut, Y. Fong and W. F. Ke, ``Composition Diamond Lemma for
 Associative Conformal Algebras", {\it J. Algebra},  \textbf{272}(2004), p.739.

\bibitem{b06}L. A. Bokut and K. P. Shum,   ``Relative Gr\"{o}bner-Shirshov Bases for Algebras
and Groups", {\it Algebra and Analisis}, \textbf{19}(2007), p.1. (in
Russian)

\bibitem{bu65}B. Buchberger, ``An Algorithm for Finding a Basis for the
Residue Class Ring of a Zero-dimensional Polynomial Ideal", Ph.D.
thesis, University of Innsbruck, Austria, (1965). (in German)

\bibitem{bu70}B. Buchberger, ``An Algorithmical Criteria for the
Solvability of Algebraic Systems of Equations", {\it Aequationes
Math.}, \textbf{4}(1970), p.374. (in German)


\bibitem{ch}E. S. Chibrikov, ``On Free Lie Conformal Algebras", {\it Vestnik Novosibirsk
State University}, \textbf{4}(2004), p.65.


\bibitem{cyz}Yuqun Chen, Yongshan Chen and Chanyan Zhong, ``Composition-Diamond Lemma for Modules",
  {\it Czech. J. Math.}, to appear. arXiv.org/abs/0804.0917.

\bibitem{ffg93} D. R. Farkas, C. D. Feustel and E. L. Green, ``Synergy
in the Theories of Gr\"{o}bner Bases and Path Algebras", {\it Can.
J. Math.}, \textbf{45}(1993), p.727.

\bibitem{H64}H. Hironaka,  ``Resolution of Singulatities of an Algebtaic Variety
over a Field if Characteristic Zero, I, II", {\it Ann. Math.},
\textbf{79}(1964), p.109, p.205.


\bibitem{H}L. Hellstr\"{o}m, ``The Diamond Lemma for Power Series Algebras",
Ph.D. thesis, 2002.

\bibitem{kl1}S.-J. Kang and K.-H. Lee, ``Gr\"{o}bner-Shirshov
Bases for Representation Theory", {\it J. Korean Math. Soc.},
\textbf{37}(2000), p.55.

\bibitem{kl3}S.-J. Kang and K.-H. Lee, ``Linear Algebraic Approach to Gr\"{o}bner-Shirshov
Bases Theory", {\it J. Algebra}, \textbf{313}(2007), p.988.

\bibitem{kobayashi} Y. Kobayashi, ``Gr\"{o}bner Bases on Algebras
Based on Well-ordered Semigroups", {\it Math. Appl. Sci. Tech.}, to
appear.

\bibitem{MikhPetr}
A.~A.~Mikhalev, ``Shirshov's Composition Techniques in Lie
Superalgebras (Non-commutative Gr\"obner Bases)", {\it Trudy Sem.
Petrovsk}, \textbf{18}(1995), p.277 (in Russian).
 English translation: {\it
J.~Math. Sci.}, \textbf{80}(1996), p.2153.

\bibitem{MV}A. A. Mikhalev and E. A. Vasilieva, ``Standard Bases of Ideals of Free
Supercommutative Polynomail Algebra ($\varepsilon$-Grobner Bases)",
{\it Proc. Second International Taiwan-Moscow Algebra Workshop},
Springer-Verlag, 2003.

\bibitem{MZ}A. A. Mikhalev and A. A. Zolotykh, ``Standard Gr\"obner-Shirshov Bases of Free Algebras over Rings, I.
Free Associative Algebras", {\it International Journal of Algebra
and Computation}, \textbf{8}(1998), p.689.

\bibitem{newman} M. H. A. Newman, ``On Theories with a Combinatorial
Definition of `Equivalence'", {\it Ann. of Math.},
\textbf{43}(1942), p.223.

\bibitem{s58}A. I. Shirshov, ``On Free Lie Rings", {\it Mat. Sb.}, \textbf{45(87)}(1958),
p.113. (in Russian)


\bibitem{s62'}A. I. Shirshov, ``Some Algorithmic Problem for $\varepsilon$-Algebras",
 {\it Sibirsk. Mat. Z.}, \textbf{3}(1962), p.132. (in Russian)

\bibitem{s62}A. I. Shirshov, ``Some Algorithmic Problem for Lie Algebras",
 {\it Sibirsk. Mat. Z.}, \textbf{3}(1962), p.292 (in Russian). English translation: {\it SIGSAM Bull.}, \textbf{33(2)}(1999), p.3.


\bibitem{st90}T. Stokes, ``Gr\"obner-Shirshov Bases in Exterior
Algebras", {\it J. Automated Reasoing}, \textbf{6}(1990), p.233.

\bibitem{Z}A. I. Zhukov, ``Complete Systems of Defining Relations in Non-associative Algebras",
 {\it Mat. Sbornik}, \textbf{69(27)}(1950), p.267.

\end{thebibliography}
\end{document}